\documentclass[12pt]{article}
\usepackage{amsmath, amssymb}

\def\cala{{\cal A}}
\def\3{\subset }
\def\4{\subseteq }
\def\<{\left<}
\def\>{\right>}
\def\vsp{\vspace*{1,5mm}\\ }
\def\n{\noindent }
\def\bk{\bigskip }
\def\bit{\begin{itemize}}
\def\eit{\end{itemize}}
\def\3{\subset }
\def\4{\subseteq }
\def\ov{\overline}

\def\0{\leqno}
\def\a{{\alpha}}

\def\barr{\begin{array}}
\def\earr{\end{array}}
\def\dd{\displaystyle}

\def\calg{{\cal G}}
\def\cals{{\cal S}}

\def\frax{\dd\frac}

\def\n{\noindent }
\def\ov#1{\overline{#1}}

\def\acom{\vrule height 1.4ex width .1ex depth -.1ex}
\def\8{\hbox{$\acom\!\raise1pt\hbox{$\times$}$}}

\title{\bf Normality degrees of finite groups}
\author{Marius T\u arn\u auceanu}
\date{December 5, 2013}

\begin{document}

\maketitle

\begin{abstract}
In this  paper we introduce and study the concept of normality
degree of a finite group $G$. This quantity measures the
probability of a random subgroup of $G$ to be normal. Explicit
formulas are obtained for some particular classes of finite
groups. Several limits of normality degrees are also computed.
\end{abstract}

\noindent{\bf MSC (2010):} Primary 20D60, 20P05;  Secondary 20D30,
20F16, 20F18.

\noindent{\bf Key words:} normality degree, subgroup lattice,
normal subgroup lattice, fixed point.

\section{Introduction}

In the last years there has been a growing interest in the use of
probability in finite group theory. One of the most important
aspects which have been studied is the probability that two
elements of a finite group $G$ commute. It is called the {\it
commutativity degree} of $G$, and has been investigated in many
papers, as \cite{3}, \cite{5}--\cite{9} or \cite{11}. Inspired by
this concept, in \cite{18} we introduced a similar notion for the
subgroups of $G$, called the {\it subgroup commutativity degree}
of $G$. This quantity is defined by $$\barr{lcl}
sd(G)&=&\frax1{|L(G)|^2}\,\left|\{(H,K)\in L(G)^2\mid
HK=KH\}\right|=\vsp &=&\frax1{|L(G)|^2}\,\left|\{(H,K)\in
L(G)^2\mid HK\in L(G)\}\right|\earr$$ (where $L(G)$ denotes the
subgroup lattice of $G$) and it measures the probability that two
subgroups of $G$ commute, or equivalently the probability that the
product of two subgroups of $G$ be a subgroup of $G$ (recall also
the natural generalization of $sd(G)$, namely the {\it relative
subgroup commutativity degree} of a subgroup of $G$, introduced
and studied in \cite{20}).

A remarkable modular sublattice of $L(G)$ is the normal subgroup
lattice $N(G)$, which consists of all normal subgroups of $G$.
Note that for an arbitrary finite group $G$ computing the number
of subgroups, as well as the number of normal subgroups, is a
difficult work. These numbers are in general unknown, excepting
for few particular classes of finite groups.

In the following we introduce the quantity
$$ndeg(G)=\dd\frac{\mid N(G) \mid}{\mid L(G) \mid}\hspace{1mm},$$ which will
be called the {\it normality degree} of $G$. Clearly, it
constitutes a signifi\-cant probabilistic aspect on subgroup
lattices of finite groups, by measuring the probability of a
random subgroup of such a group to be normal. The normality degree
is closely connected to a special type of an action of a finite
group on a lattice, introduced and studied in \cite{15}. Recall
that, given a finite group $G$ acting on a lattice $(L, \wedge,
\vee)$, we say that $L$ is a $G$-$lattice$ if the following two
equalities hold
$$\begin{array}{l}
g \circ(l \wedge l')=(g \circ l) \wedge (g \circ l'),\\
g \circ(l \vee l')=(g \circ l) \vee (g \circ l'),\end{array}$$ for
all $g\in G$ and $l, l' \in L$, that is the action $\circ$ of $G$
on $L$ is compatible with the binary operations $\wedge$ and
$\vee$ of $L$. For a finite $G$-lattice $L$, the set $Fix_G(L)=\{l
\in L \mid g \circ l=l,\hspace{1mm} {\rm for} \hspace{1mm} {\rm
all} \hspace{1mm} g \in G\}$ forms a $G$-sublattice of $L$ and the
quantity $$\dd\frac{\mid Fix_G(L) \mid}{\mid L \mid}\0(*)$$
measures the probability of an element of $L$ to be fixed with
respect to $\circ$. Moreover, if we assume that both the initial
element and the final element of $L$ are contained in $Fix_G(L)$,
then the map $f_L : L \longrightarrow L$ defined by
$f_L(l)=\dd\bigwedge_{g \in G}g \circ l$, for any $l \in L$, is
isotone. Therefore the set $Fix(f_L)=\{l \in L \mid f_L(l)=l\}$ is
also a $G$-sublattice of $L$, according to the well-known
fixed-point theorem of complete lattices. Again, a specific
quantity associated to $L$, namely
$$\dd\frac{\mid Fix(f_L) \mid}{\mid L \mid}\hspace{1mm},\0(**)$$ can be seen as
a probabilistic aspect on $L$, more precisely it measures the
probability of an element of $L$ to be a fixed point of $f_L$. One
of the most important examples of a $G$-lattice is constituted by
the subgroup lattice $L(G)$ associated to $G$. In this case the
action of $G$ on $L(G)$ is defined by $g \circ H=H^g$, for all
$(g, H) \in G \times L(G)$, and $f_{L(G)}$ maps every subgroup $H
\in L(G)$ into its core in $G$. Then both $G$-sublattices
$Fix_G(L(G))$ and $Fix(f_{L(G)})$ of $L(G)$ will coincide with the
normal subgroup lattice $N(G)$. In other words, both quantities
($*$) and ($**$) are equal to the normality degree $ndeg(G)$ of
$G$. Hence $ndeg(G)$ measures the probability of a random subgroup
of $G$ to be a fixed point of $L(G)$ relative to the above
canonical action of $G$ on $L(G)$, and also to be a fixed point of
the map $f_{L(G)}$.

All the previous remarks give a strong motivation to study the
normality degree of finite groups. In our paper a first step of
this study is made.

The paper is organized as follows. Some basic properties and
results on normality degree are presented in Section 2. Section 3
deals with normality degrees for two special classes of finite
groups: semidirect products of finite cyclic groups and finite
$p$-groups possessing a cyclic maximal subgroup. An interesting
density result of normality degree is proved in Section 4. In the
final section several conclusions and further research directions
are indicated.

Most of our notation is standard and will usually not be repeated
here. Elementary notions and results on lattices (respectively on
groups) can be found in \cite{2} (respectively in \cite{4} and
\cite{14}). For subgroup lattice concepts we refer the reader to
\cite{12}, \cite{15} and \cite{16}.

\section{Basic properties of normality degree}

Let $G$ be a finite group. First of all, remark that the normality
degree $ndeg(G)$ satisfies the following relation
$$0<ndeg(G)\le 1.$$ Moreover, we have $ndeg(G)=1$ if and only if all subgroups of $G$ are normal,
that is $G$ is a Dedekind group. As we have seen in \cite{18}, the
normality degree and the subgroup commutativity degree of $G$ are
connected by the inequality $$ndeg(G)\le sd(G).\0(1)$$ Clearly,
this becomes an equality if and only if for every subgroup $H$ of
$G$ the set $C(H)$ consisting of all subgroups of $G$ which
commute with $H$ coincides with $N(G)$. Since $H$ itself is
contained in $C(H)$, it must be normal and so $G$ is a Dedekind
group. Hence the following result holds.

\bk\n{\bf Proposition 2.1.} {\it Let $G$ be a finite group. Then
the following conditions are equivalent: \bit
\item[\rm a)]$ndeg(G)=1$.
\item[\rm b)] $ndeg(G)=sd(G).$
\item[\rm c)] $G$ is a Dedekind group.\eit}
\bigskip

Next, let $\cals$ be a set of representatives for the conjugacy
classes of subgroups of $G$ with at least two elements. Then
$$|L(G)|=|N(G)|+\dd\sum_{H \in \cals}(G:N_G(H)),$$ which implies
that $$ndeg(G)=\dd\frax{|N(G)|}{|N(G)|+\dd\sum_{H \in
\cals}(G:N_G(H))}\hspace{1mm}.\0(2)$$ This equality can be used to
calculate the normality degree of finite groups whose conjugacy
classes of subgroups are completely determined. The simplest
examples are constituted by the symmetric groups $S_3$ and $S_4$,
for which one obtains
$$ndeg(S_3)=\dd\frax{1}{2} \hspace{2mm} {\rm and} \hspace{1mm} ndeg(S_4)=\dd\frax{2}{15} \hspace{1mm}.$$
In particular, if $G$ is a finite $p$-group, (2) leads us to an
inequality satisfied by $ndeg(G)$, namely
$$ndeg(G)\le \dd\frax{|N(G)|}{|N(G)|+p\hspace{1mm}|\cals|}\hspace{1mm}.$$

In many situations computing the normality degree of a finite
group is reduced to computing the number of all its subgroups. One
of them is constituted by finite groups with few normal subgroups,
as the symmetric groups.

\bk\n{\bf Example 2.2.} The following equality holds
$$ndeg(S_n)=\dd\frax{3}{|L(S_n)|}\hspace{1mm}, \hspace{1mm}{\rm
for}\hspace{1mm} {\rm all}\hspace{1mm} n \geq 5.$$ Mention that we
also have $$ndeg(S_n \times S_n)=\dd\frax{10}{|L(S_n \times
S_n)|}\hspace{1mm}, \hspace{1mm}{\rm for}\hspace{1mm} {\rm
all}\hspace{1mm} n \geq 5,$$ and a formula for $ndeg(S_n^k)$ which
depends only on $|L(S_n^k)|$ can be easily inferred, according to
\cite{10}.\bk

In the following assume that $G$ and $G'$ are two finite groups.
It is obvious that if $G\cong G'$, then $ndeg(G)=ndeg(G')$. In
particular, we infer that any two conjugate subgroups of a finite
group have the same normality degree. The above conclusion is not
true in the case when $G$ and $G'$ are only lattice-isomorphic, as
shows the next example.

\bk\n{\bf Example 2.3.} Let $G$ be the finite elementary abelian
$3$-group $\mathbb{Z}^n_3$ (where $n\ge2$) and $G'$ be a
semidirect product of an elementary abelian normal subgroup $A$ of
order $3^{n-1}$ by the group $B \cong \mathbb{Z}_2$ which induces
a nontrivial power automorphism on $A$. Then both $G$ and $G'$ are
contained in the class $P(n,3)$ (see \cite{12}, page 49) and so
they are $L$-isomorphic. We have $ndeg(G)=1$, because $G$ is
abelian. On the other hand, in Section 2 of \cite{18} we have
proved that $sd(G')<1$. This implies that $ndeg(G')<1$, in view of
(1). Hence $ndeg(G)\ne ndeg(G').$ \bk

By a direct calculation we obtain
$$ndeg(S_3 \times \mathbb{Z}_2)=\dd\frax{7}{16} \ne \frax{1}{2}=ndeg(S_3)ndeg(\mathbb{Z}_2)$$
and therefore in general we don't have $ndeg(G\times
G')=ndeg(G)ndeg(G')$. Clearly, a sufficient condition in order to
this equality holds is $$(\mid G \mid,\mid G'\mid)=1,$$ that is
$G$ and $G'$ are of coprime orders. This remark can be extended to
an arbitrary finite direct product.

\bk\n{\bf Proposition 2.4.} {\it Let $(G_i)_{i=\overline{1,k}}$ be
a family of finite groups having coprime orders. Then
$$ndeg(\prod_{i=1}^k G_i)=\prod_{i=1}^k ndeg(G_i).$$}

The following immediate consequence of Proposition 2.4 shows that
computing the normality degree of a finite nilpotent group is
reduced to finite $p$-groups.

\bk\n{\bf Corollary 2.5.} {\it If $G$ is a finite nilpotent group
and $(G_i)_{i=\ov{1,k}}$ are the Sylow subgroups of $G$, then
$$ndeg(G)=\prod_{i=1}^k ndeg(G_i).$$}

\section{Normality degrees for some\\ classes of finite groups}

In this section we determine explicitly the normality degree of
several finite groups. The most significant results are obtained
for the class of finite dihedral groups and for the class of
finite $p$-groups possessing a cyclic maximal subgroup.

\subsection{The normality degree of some semidirect\\ products of finite groups}

Let $p$ be a prime, $n \ge 2$ be an integer such that $p \nmid n$
and $f:\mathbb{Z}_p \longrightarrow{\rm Aut}(\mathbb{Z}_n)$ be a
group homomorphism. Put $\hat k_0=f(\bar 1)(\hat1)$ and suppose
that $k_0 \neq 1$. Then we have $(k_0,n)=1$ and
$$f(\bar x)(\hat y)=k_0^x\hat y,\ \mbox{for any $\bar x \in \mathbb{Z}_p,\ \hat y \in \mathbb{Z}_n.$}$$
Denote by $G$ be the semidirect product of $\mathbb{Z}_p$ and
$\mathbb{Z}_n$ with respect to $f$. Recall that the operation of
$G$ is defined by
$$(\bar x_1,\hat y_1)\cdot(\bar x_2,\hat y_2)=(\bar x_1+\bar
x_2,k^{x_2}_0\hat y_1+\hat y_2), \mbox{ for all $(\bar x_1,\hat
y_1),(\bar x_2,\hat y_2)\in G.$}$$ It is well-known that the maps
$$\begin{array}{lll}
\sigma_1:\mathbb{Z}_p \longrightarrow G,& \sigma_1(\bar x)=(\bar
x,\hat0),& \mbox{ for any }\bar x \in \mathbb{Z}_p,\vsp
\sigma_2:\mathbb{Z}_n \longrightarrow G,& \sigma_2(\hat y)=(\bar
0,\hat y),& \mbox{ for any }\hat y \in \mathbb{Z}_n,\end{array}$$
are injective group homomorphisms. Moreover, if
$H=\sigma_1(\mathbb{Z}_p)$ and $K=\sigma_2(\mathbb{Z}_n),$ then
$H$ is a subgroup of $G$ and $K$ is a normal subgroup of $G$,
which satisfy
$$G=HK,\ H \cap K=\{(\bar0,\hat0)\}.$$

In the following our goal is to compute explicitly the normality
degree of $G$. First of all, we shall give a precise description
of $L(G)$ (for more details, see Section 3.2 of \cite{15}). Let
$G_1$ be a subgroup of $G$. Then $\mid G_1 \mid$ is a divisor of
$pn.$

In the case when $p \nmid \hspace{2mm} \mid G_1 \mid$ we shall
prove that $G_1 \subseteq K.$ Indeed, if we assume that $G_1
\nsubseteq K,$ then we have $K \subset G_1K \subseteq G$ and so
the index $(G_1K:K)$ of $K$ in $G_1K$ is $\ge 2$. Since
$p=(G:K)=(G:G_1K)(G_1K:K)$ is prime, one obtains $(G_1K:K)=p$ and
$(G:G_1K)=1,$ i.e. $G_1K=G$. It results
$$ G_1/ G_1 \cap K \cong G_1K/K = G/K,$$
which shows that $\mid G_1/ G_1 \cap K \mid =p$\ and therefore $p
\mid \hspace{2mm}\mid G_1 \mid,$ a contradiction. Hence
$$G_1 \in L(K)=L(\sigma_2(\mathbb{Z}_n))=\sigma_2(L(\mathbb{Z}_n)).\0(3)$$

In the case when $p \mid \hspace{2mm}\mid G_1 \mid$ at least a
subgroup of order $p$ is contained in $G_1$. Let $\{H_1=H,
H_2,...,H_{n_p}\}$ be the set of all Sylow $p$-subgroups of $G$,
where $n_p=(G:N_G(H))$. By a direct calculation, the normalizer
$N_G(H)$ of $H$ in $G$ can be easily determined.

\bk\n{\bf Lemma 3.1.1.} {\it The following equality holds
$$N_G(H)=\{\hspace{1mm}(\bar x,\hat y) \in G \mid \bar
x \in \mathbb{Z}_p,\ \hat y\!\in \langle\dd\frac{\hat
n}d\rangle\},$$ where $d =(k_0-1,n).$}\bk

Then $n_p=\dd\frac{n}{d}=\dd\frac{n}{(k_0-1,n)}\hspace{1mm}.$ For
every $i \in \{1,2,...,n_p\}$ there exists $z_i \in G$
$(z_1=(\bar0,\hat0))$ such that $H_i=H^{z_i}.$ One obtains
$$G=G^{z_i}=(HK)^{z_i}=H^{z_i}K^{z_i}=H_iK,\hspace{1mm} i=\ov{1,n_p}\hspace{1mm}.$$
Suppose that $H_i \subseteq G_1$ for some $i \in \{1,2,...,n_p\}$.
It results $G_1=G_1 \cap G=G_1 \cap (H_iK)=H_i(G_1 \cap K)$ and
thus $G_1$ is contained in the set
$$\cala\hspace{-1mm}=\hspace{-1mm}\{H^{z_i} \sigma_2 (\langle\dd\frac{\hat
n}{k}\rangle)\hspace{-1mm}\mid k \hspace{-1mm}\mid \hspace{-1mm}n,
\hspace{1mm}i=\ov{1,n_p}\hspace{1mm}\}\hspace{-1mm}=\hspace{-1mm}\{(H
\sigma_2 (\langle\dd\frac{\hat
n}{k}\rangle))^{z_i}\hspace{-1mm}\mid  k \hspace{-1mm}\mid
\hspace{-1mm}n, \hspace{1mm}i=\ov{1,n_p}\hspace{1mm}\}.\0(4)$$

In order to determine the number of elements of $\cala$, we need
to compute the normalizer in $G$ of such an element.

\bk\n{\bf Lemma 3.1.2.} {\it If $k$ is a divisor of $n$, then
$$N_G(H \sigma_2 (\langle\dd\frac{\hat n}{k}\rangle))=
\{\hspace{1mm}(\bar x,\hat y) \in G \mid \bar x\in\mathbb{Z}_p,\
\hat y \in \langle\dd\frac{\hat n}{\varepsilon(k)}\rangle\},$$
where $\varepsilon(k)=(k(k_0-1),n).$}\bk

From Lemma 3.1.2 we easily infer that
$$\mid \cala \mid=\dd\sum_{k \mid n} \dd\frac{\varepsilon(k)}{(k_0-1,n)}=\dd\sum_{k \mid n} \hspace{1mm}(k,\dd\frac{n}{(k_0-1,n)}).\0(5)$$

Now, by using the relations (3), (4) and (5), we are able to
describe the subgroup structure of $G$.

\bk\n{\bf Proposition 3.1.3.} {\it The subgroup lattice $L(G)$ of
the above semidirect pro-duct $G$ is given by the following
equality:
$$L(G)=\sigma_2(L(\mathbb{Z}_n)) \cup \cala.$$
Moreover, the total number of subgroups of $G$ is
$$\mid L(G) \mid=\tau(n)+\dd\sum_{k \mid n} \hspace{1mm}(k,\dd\frac{n}{(k_0-1,n)}),$$
where $\tau(n)$ denotes the number of all divisors of $n$.}\bk

Clearly, the normal subgroups of $G$ are all subgroups contained
in $K$ and $G$ itself, that is
$$N(G)=\sigma_2(L(\mathbb{Z}_n)) \cup \{G\},$$ and therefore
$$\mid N(G) \mid=\tau(n)+1.$$

Hence we have proved the following theorem.

\bk\n{\bf Theorem 3.1.4.} {\it The normality degree of the above
semidirect product $G$ is given by the following equality:
$$ndeg(G)=\dd\frac{\tau(n)+1}{\tau(n)+\dd\sum_{k \mid n} \hspace{1mm}(k,\dd\frac{n}{(k_0-1,n)})}\hspace{1mm}.\0(6)$$}

\n{\bf Remark.} Let $r=\dd\frac{n}{(k_0-1,n)}$. Then $1 \leq (k,r)
\leq k,r$, for all divisors $k$ of $n$. So, by (6) we infer that
$ndeg(G)$ satisfies the following inequalities:
$$\hspace{-20mm}ndeg(G) \leq \dd\frac{\tau(n)+1}{2 \tau(n)}\hspace{1mm},\0(7)$$
$$\hspace{-14mm}ndeg(G) \geq \dd\frac{\tau(n)+1}{\tau(n)+\sigma(n)}\hspace{1mm},\0(8)$$
$$ndeg(G) \geq \dd\frac{\tau(n)+1}{\tau(n)(r+1)} > \dd\frac{1}{r+1}\hspace{1mm}.\0(9)$$
\smallskip

Next, let us assume that $p=2$ and $k_0=n-1$. Then the group $G$
studied above is the dihedral group $D_{2n}$. Recall that $D_{2n}$
is the symmetry group of a regular polygon with $n$ sides and it
has the order $2n$. The most convenient abstract description of
$D_{2n}$ is obtained by using its generators: a rotation $x$ of
order $n$ and a reflection $y$ of order $2$. Under these
notations, we have
$$D_{2n}=\langle x,y\mid x^n=y^2=1,\ yxy=x^{-1}\rangle.$$

Since $n$ is odd, it results $(k_0-1,n)=(n-2,n)=1$ and so
$$\dd\sum_{k \mid n} \hspace{1mm}(k,\dd\frac{n}{(k_0-1,n)})=\dd\sum_{k \mid n}
\hspace{1mm}(k,n)=\sigma(n),$$ where $\sigma(n)$ denotes the sum
of all divisors of $n$. Thus, (6) leads us to
$$ndeg(D_{2n})=\dd\frac{\tau(n)+1}{\tau(n)+\sigma(n)}\hspace{1mm},\0(10)$$ that
is (8) becomes an equality for $G=D_{2n}$ with $n$ odd.\bk

A similar formula can be also obtained for even positive integers
$n$. In this case it is well-known that we have
$$N(D_{2n})=L(\langle x\rangle) \cup \{D_{2n}, \langle x^2,y\rangle, \langle x^2,xy\rangle\}$$
and therefore
$$ndeg(D_{2n})=\dd\frac{\tau(n)+3}{\tau(n)+\sigma(n)}\hspace{1mm}.\0(11)$$

Hence (10) and (11) imply the following result.

\bk\n{\bf Corollary 3.1.5.} {\it The normality degree of the
dihedral group $D_{2n}$ is given by the following equality:
$$ndeg(D_{2n})=\left\{\barr{lll}
\dd\frac{\tau(n)+1}{\tau(n)+\sigma(n)}\hspace{1mm},& {\rm if}\hspace{1mm} n \hspace{1mm} {\rm is} \hspace{1mm} {\rm odd}\\
&&\\
\dd\frac{\tau(n)+3}{\tau(n)+\sigma(n)}\hspace{1mm},& {\rm
if}\hspace{1mm} n \hspace{1mm} {\rm is} \hspace{1mm} {\rm even
\hspace{1mm}.}\earr\right.\0(12)$$}

\n{\bf Remark.} A simple arithmetic exercise shows that $\tau(n)+2
\leq \sigma(n)$, for all odd positive integers $n \neq 1$, and
$\tau(n)+6 \leq \sigma(n)$, for all even positive integers $n \neq
2,4$. These inequalities give us an upper bound for the normality
degree of $D_{2n}$, namely
$$ndeg(D_{2n}) \leq \dd\frac{1}{2}\hspace{1mm},$$for all $n \neq
2,4.$ Mention also that we have $ndeg(D_{2n})=\dd\frac{1}{2}$ if
and only if $n=3$, that is in the class of finite dihedral groups
only $D_6$ (which is isomorphic to $S_3$) has the normality degree
$\dd\frac{1}{2}\hspace{1mm}.$ \bk

In the end of this subsection, we note that the normality degrees
of other semidirect products of finite groups can be also
computed. Such an example is constituted by ZM-groups, that is
finite groups with all Sylow subgroups cyclic. It is well-known
(see \cite{4}, I) that a ZM-group possesses a presentation of type
$${\rm ZM}(m,n,r)=\langle a, b \mid a^m = b^n = 1, \hspace{1mm}b^{-1} a
b = a^r\rangle\hspace{1mm},$$where the triple $(m,n,r)$ satisfies
the conditions
$$(m,n) = (m, r-1) = 1 \mbox{ and }
r^n \equiv 1 \hspace{1mm}({\rm mod}\hspace{1mm}m).$$The subgroups
of ${\rm ZM}(m,n,r)$ have been computed in \cite{1}: $$\mid L({\rm
ZM}(m,n,r))\mid=\dd\sum_{m_1\mid m}\sum_{n_1\mid
n}\hspace{1mm}(m_1,\dd\frac{r^n-1}{r^{n_1}-1})\hspace{1mm},$$while
the number of normal subgroups of ${\rm ZM}(m,n,r)$ has been
determined in \cite{21}:
$$\mid N({\rm
ZM}(m,n,r))\mid=\dd\sum_{n_1\mid
n}\tau((m,r^{n_1}-1))\hspace{1mm}.$$In this way, one obtains
$$ndeg({\rm ZM}(m,n,r))=\dd\frac{\dd\sum_{n_1\mid n}\tau((m,r^{n_1}-1))}{\dd\sum_{m_1\mid m}\sum_{n_1\mid n}\hspace{1mm}(m_1,\dd\frac{r^n-1}{r^{n_1}-1})}\hspace{1mm}.\0(13)$$Finally,
remark that, by taking $n=2$ and $r=m-1$ with $m$ odd in (13), the
previous formula (10) is obtained. This is not a surprise,
according to the group isomorphism ${\rm ZM}(m,2,m-1)\cong
D_{2m}$.

\subsection{The normality degree of finite $p$-groups\\ possessing a cyclic maximal subgroup}

Let $p$ be a prime, $n \ge 3$ be an integer and denote by $\calg$
the class consisting of all finite $p$-groups of order $p^n$
having a maximal subgroup which is cyclic. Obviously, $\calg$
contains finite abelian $p$-groups of type $\mathbb{Z}_p \times
\mathbb{Z}_{p^{n-1}}$ whose normality degree is 1, but some finite
nonabelian $p$-groups belong to $\calg$, too. They are
exhaustively described in Theorem 4.1, \cite{14}, II: a nonabelian
group is contained in $\calg$ if and only if it is isomorphic to
\bit\item[--] $M(p^n)=\langle x,y\mid x^{p^{n-1}}=y^p=1,\ y^{-1}x
y=x^{p^{n-2}+1}\rangle$\eit when $p$ is odd, or to one of the
following groups \bit\item[--] $M(2^n)\ (n \ge 4),$
\item[--] the dihedral group $D_{2^n}$,
\item[--] the generalized quaternion group
$$Q_{2^n}=\langle x,y\mid x^{2^{n-1}}=y^4=1,\ yxy^{-1}=x^{2^{n-1}-1}\rangle,$$
\item[--] the quasi-dihedral group
$$S_{2^n}=\langle x,y\mid x^{2^{n-1}}=y^2=1,\ y^{-1}xy=x^{2^{n-2}-1}\rangle\ (n \ge 4)$$\eit
when $p=2$. \bk

In the following the normality degrees of the above $p$-groups
will be determined. We recall first the explicit formulas for the
total number of subgroups of these groups, found in \cite{18}.

\bk\n{\bf Lemma 3.2.1.} {\it The following equalities hold: \bit
\item[\rm a)] $\mid L(M(p^n)) \mid \hspace{1mm}=(1+p)n+1-p \hspace{1mm},$
\item[\rm b)] $\mid L(D_{2^n}) \mid \hspace{1mm}=2^n+n-1 \hspace{1mm},$
\item[\rm c)] $\mid L(Q_{2^n}) \mid \hspace{1mm}=2^{n-1}+n-1 \hspace{1mm},$
\item[\rm d)] $\mid L(S_{2^n}) \mid \hspace{1mm}=3 \cdot 2^{n-2}+n-1 \hspace{1mm}.$
\eit}

In order to compute the normality degree of the nonabelian
$p$-groups that belong to $\calg$, we need to know the number of
their normal subgroups. Our reasoning is founded on the following
simple remark: such a group $G$ possesses a unique normal subgroup
of order $p$, namely $\langle x^q\rangle$ (where $q=p^{n-2}$ and
$x$ denotes a generator of a cyclic maximal subgroup of $G$). We
infer that there exists a bijection between the set of nontrivial
normal subgroups of $G$ and $N(G/\langle x^q\rangle)$, that is
$$\mid N(G) \mid\hspace{1mm}=1+ \mid N(G/\langle x^q\rangle) \mid.\0(14)$$

For $G=M(p^n)$, the minimal normal subgroup $\langle x^q\rangle$
is in fact the commutator subgroup $D(M(p^n))$ of $M(p^n)$ and we
have
$$M(p^n)/D(M(p^n)) \cong \mathbb{Z}_p \times
\mathbb{Z}_{p^{n-2}}\hspace{1mm}.$$ Since $\mathbb{Z}_p \times
\mathbb{Z}_{p^{n-2}}$ is abelian, the number of its normal
subgroups coincides with the number of all its subgroups. Put
$x_n=\mid L(\mathbb{Z}_p \times
\mathbb{Z}_{p^{n-2}})\mid\hspace{-1mm}.$ This number can be easily
determined by using the following lemma, established in \cite{17}
(see also \cite{19}).

\bk\n{\bf Lemma 3.2.2.} {\it  For every $0\le\a\le\a_1+\a_2$,  the
number of all subgroups of order $p^{\a_1+\a_2-\a}$ in the finite
abelian $p$-group $\mathbb{Z}_{p^{\a_1}} \times
\mathbb{Z}_{p^{\a_2}}\hspace{1mm} (\a_{1}\leq\a_{2})$ is
$$\left\{\barr{lll}
\displaystyle\frac{p^{\a+1}-1}{p-1}\,,&\mbox{ if }&0\le\a\le\a_1\vspace*{3mm}\\
\displaystyle\frac{p^{\a_1+1}-1}{p-1}\,,&\mbox{ if }&\a_1\le\a\le\a_2\vspace*{3mm}\\
\displaystyle\frac{p^{\a_1+\a_2-\a+1}-1}{p-1}\,,&\mbox{ if
}&\a_2\le\a\le\a_1+\a_2.\earr\right.$$\bigskip In particular, the
total number of subgroups of $\mathbb{Z}_{p^{\a_1}} \times
\mathbb{Z}_{p^{\a_2}}$ is
$$\dd\frac{1}{(p{-}1)^2}\left[(\a_2{-}\a_1{+}1)p^{\a_1{+}2}{-}(\a_2{-}\a_1{-}1)p^{\a_1{+}1}{-}(\a_1{+}\a_2{+}3)p{+}(\a_1{+}\a_2+1)\right].$$}\\
By taking $\a_1=1$ and $\a_2=n-2$, one obtains
$$x_n=\dd\frac{1}{(p{-}1)^2}\left[(n-2)p^3{-}(n-4)p^2{-}(n+2)p{+}n\right]=(1+p)n-2p.$$ So, (14) becomes
$$\mid N(M(p^n)) \mid\hspace{1mm}=1+x_n\hspace{1mm}=(1+p)n+1-2p.$$

For every $G \in \{D_{2^n}, Q_{2^n}, S_{2^n}\}$ the minimal normal
subgroup $\langle x^q\rangle$ coincides with the center $Z(G)$ of
$G$ and we have
$$G/Z(G) \cong D_{2^{n-1}},$$ therefore $$\mid N(G) \mid\hspace{1mm}=1+ \mid N(D_{2^{n-1}})
\mid.$$ Let $G=D_{2^{n-1}}$ in the above equality and set $y_n=
\mid N(D_{2^{n-1}}) \mid\hspace{-1mm}.$ Then the integer sequence
$(y_n)_{n \in \mathbb{N^*}}$ satisfies the recurrence relation
$y_n=1+y_{n-1},$ which shows that $y_n=n+3, \hspace{1mm} {\rm for}
\hspace{1mm} {\rm any} \hspace{1mm} n \in \mathbb{N^*}.$ Thus
$$\mid N(G) \mid\hspace{1mm}=1+y_{n-1}=y_n=n+3,$$ for all above
2-groups $G$. Hence we have proved the following result.

\bk\n{\bf Lemma 3.2.3.} {\it The following equalities hold: \bit
\item[\rm a)] $\mid N(M(p^n)) \mid\hspace{1mm} =(1+p)n+1-2p \hspace{1mm},$
\item[\rm b)] $\mid N(G) \mid\hspace{1mm}=n+3$, for all $G \in \{D_{2^n}, Q_{2^n}, S_{2^n}\}.$
\eit}
\bigskip

Now, it is clear that Lemma 3.2.1 and Lemma 3.2.3 imply the next
theorem.

\bk\n{\bf Theorem 3.2.4.} {\it The normality degrees of the
nonabelian groups in the class $\calg$ are given by the following
equalities: \bit
\item[\rm 1)] $ndeg(M(p^n))=\frax{(1+p)n+1-2p}{(1+p)n+1-p}\hspace{1mm},$
\item[\rm 2)] $ndeg(D_{2^n})=\frax{n+3}{2^n+n-1}\hspace{1mm},$
\item[\rm 3)] $ndeg(Q_{2^n})=\frax{n+3}{2^{n-1}+n-1}\hspace{1mm},$
\item[\rm 4)] $ndeg(S_{2^n})=\frax{n+3}{3 \cdot 2^{n-2}+n-1}\hspace{1mm}.$
\eit}
\bigskip

\bk\n{\bf Remark.} The normality degree of the dihedral group
$D_{2^n}$ can also be directly computed by using Corollary 3.1.5:
$$ndeg(D_{2^n})=\dd\frac{\tau(2^{n-1})+3}{\tau(2^{n-1})+\sigma(2^{n-1})}=
\frax{n+3}{2^n+n-1}\hspace{1mm}.$$
\bigskip

The following limits are immediate from Theorem 3.2.4.

\bk\n{\bf Corollary 3.2.5.} {\it We have: \bit
\item[\rm a)]$\dd\lim_{n\to\infty}ndeg(M(p^n))=1$, for any fixed prime $p$.
\item[\rm b)]$\dd\lim_{n\to\infty}ndeg(G)=0$, for all $G \in \{D_{2^n}, Q_{2^n},
S_{2^n}\}.$\eit}\bk

We end this section by mentioning that the normality degree of any
finite nilpotent group whose Sylow subgroups belong to $\calg$ can
explicitly be calculated, in view of Corollary 2.5.

\section{A density result of normality degree}

As we have seen in Section 3, there are some sequences of finite
groups $(G_n)_{n \in \mathbb{N}}$ satisfying
$\dd\lim_{n\to\infty}ndeg(G_n) \in \{0,1\}.$ In this section our
purpose is to extend this result by proving that each $x$ in the
interval [0,1] is the limit of a certain sequence of normality
degrees of finite groups.\bk

First of all, we shall prove the above result for rational numbers
in [0,1].

\bk\n{\bf Theorem 4.1.} {\it For every $x \in [0,1] \cap
\mathbb{Q}$ there exists a sequence $(G_n)_{n \in \mathbb{N}}$ of
finite groups such that $\dd\lim_{n\to\infty}ndeg(G_n)=x.$}\bk

\n{\bf Proof.} For $x=0$ and $x=1$ our statement is already
verified in the previous section, by taking $G_n=D_{2^n}$ (or
$G_n=Q_{2^n}$, or $G_n=S_{2^n}$) and $G_n=M(p^n)$, respectively.
Let $x \in (0,1) \cap \mathbb{Q}$. Then $x=\dd\frax{a}{b}$, where
$a, b \in \mathbb{N}^*$ and $a < b$. Denote by $(p_n)_{n \in
\mathbb{N}}$ the sequence of the prime numbers and choose the
disjoint strictly increasing subsequences $(k_n^1), (k_n^2), ...,
(k_n^{b-a})$ of $\mathbb{N}$. We also consider
$G_i=M(p_{k_n^i}^{a+i+1})$, $i=1,2, ..., b-a$. Then the normality
degree of $G_i$ is given by
$$ndeg(G_i)=\dd\frax{(a+i-1)p_{k_n^i}+a+i+2}{(a+i)p_{k_n^i}+a+i+2}$$
and we have
$$\dd\lim_{n\to\infty}ndeg(G_i)=\dd\frax{a+i-1}{a+i}\hspace{1mm},$$ for all
$i=\ov{1,b-a}$. Let $G=\prod_{i=1}^{b-a} G_i$. From Corollary 2.5
it results $$ndeg(G)=\prod_{i=1}^{b-a} ndeg(G_i).$$Hence
$$\dd\lim_{n\to\infty}ndeg(G_n)=\prod_{i=1}^{b-a}\dd\lim_{n\to\infty}ndeg(G_i)=\prod_{i=1}^{b-a}\dd\frax{a+i-1}{a+i}=\dd\frax{a}{b}=x,$$
which completes our proof.\hfill\rule{1,5mm}{1,5mm}\bk

Since the set $[0,1] \cap \mathbb{Q}$ is dense in [0,1], by
Theorem 4.1 we infer the following corollary.

\bk\n{\bf Corollary 4.2.} {\it For every $x \in [0,1]$ there
exists a sequence $(G_n)_{n \in \mathbb{N}}$ of finite groups such
that $\dd\lim_{n\to\infty}ndeg(G_n)=x.$}\bk

Let $a, b \in \mathbb{N}^*$ with $a < b$. In general, there is no
finite group $G$ satisfying both equalities
$$|N(G)|=a \hspace{1mm}{\rm and}\hspace{1mm} |L(G)|=b.$$ The above
system has no solution $G$ even in the particular case when
$b=a+1$. In contrast with this statement, for several values of
$a$ we are able to construct finite groups $G$ such that
$$ndeg(G)=\dd\frax{a}{a+1}\hspace{1mm}.$$ For example, we have
$ndeg(S_3)=\dd\frax{1}{2}$ and $ndeg(M(5^4))=\dd\frax{3}{4}$ (more
generally, a fraction $\dd\frax{a}{a+1}$ is the normality degree
of a finite $p$-group of type $M(p^n)$ if and only if there is a
prime $q$ such that $q+1$ divides $a+3$). Inspired by these
examples, we came up with the following conjecture.

\bk\n{\bf Conjecture 4.3.} {\it For every $a \in \mathbb{N}^*$
there exists a finite group $G$ such that
$ndeg(G)=\dd\frax{a}{a+1}$\hspace{1mm}.}\bk

Finally, notice that it is natural to generalize Conjecture 4.3 in
the following manner.

\bk\n{\bf Conjecture 4.4.} {\it For every $x \in (0,1] \cap
\mathbb{Q}$ there exists a finite group $G$ such that
$ndeg(G)=x$.}

\section{Conclusions and further research}

All our previous results show that the concept of normality degree
introduced in this paper can constitute a significant aspect of
probabilistic finite group theory. It is clear that the study
started here can successfully be extended to other classes of
finite groups. This will surely be the subject of some further
research. \bk

Two interesting conjectures on normality degree have been
formulated in Section 4. Another open problems concerning this
topic are the following:

\bk\n{\bf Problem 4.1.} Given a finite group $G$, a subgroup $H$
of $G$ and a normal subgroup $N$ of $G$, which is the connection
between $ndeg(G)$ and $ndeg(H)$, respectively between $ndeg(G)$
and $ndeg(G/N)$?

\bk\n{\bf Problem 4.2.} Give explicit formulas for the normality
degrees of other classes of finite groups.

\bk\n{\bf Problem 4.3.} For a fixed $a \in (0,1)$, describe the
structure of finite groups $G$ satisfying
$$ndeg(G)=(\le,\hspace{1mm}\ge)\hspace{1mm}a.$$

\bk\n{\bf Problem 4.4.} Study the properties of the map $ndeg$
from the class of finite groups to [0,1]. What can be said about
two finite groups having the same normality degree?

\bk\n{\bf Problem 4.5.} As we have seen in Corollary 3.2.5, the
following equalities hold
$$\lim_{n\to\infty}ndeg(D_{2^{n-1}})=\lim_{n\to\infty}ndeg(Q_{2^{n-1}})=\lim_{n\to\infty}ndeg(S_{2^{n-1}})=0.$$
Is this true for other "natural" classes of finite groups?

\vspace*{5ex}\small

\hfill
\begin{minipage}[t]{5cm}
Marius T\u arn\u auceanu \\
Faculty of  Mathematics \\
``Al.I. Cuza'' University \\
Ia\c si, Romania \\
e-mail: {\tt tarnauc@uaic.ro}
\end{minipage}

\end{document}